\documentclass[11pt, twoside]{article}
\usepackage{amsmath}
\usepackage{amsfonts}
\usepackage{amssymb}
\usepackage{enumerate}

\usepackage{array}
\usepackage{color}

\usepackage[a4paper]{geometry}
\usepackage{hyperref}
\usepackage{fancyhdr}
\usepackage{graphicx}
\usepackage{tikz}
\usepackage{endnotes}

\let\footnote=\endnote
\begin{document}

\title{The L\'evy-It\^o Decomposition Theorem}
\author{by J.L. Bretagnolle\\Translation and notes by P. Ouwehand}\date{}

\maketitle
\begin{abstract}This a free translation with additional explanations of {\em Processus \`a Accroissement Independants Chapitre I: La D\'ecomposition de Paul L\'evy}, by J.L. Bretagnolle, in {\em Ecole d'Et\'e de Probabilit\'es}, Lecture Notes in Mathematics 307, Springer 1973. The L\'evy-Khintchine representation of infinitely divisible distributions is obtained as a by-product. 

 As this proof makes use of martingale methods, it is pedagogically more suitable for students of financial mathematics than some other approaches. It is hoped that the end notes will also help to make the proof more  accessible to this audience.\end{abstract}


\pagestyle{fancy}
\fancyhf{}

\newtheorem{theorem}{Theorem}[section]
\newtheorem{lemma}[theorem]{Lemma}
\newtheorem{corollary}[theorem]{Corollary}
\newtheorem{proposition}[theorem]{Proposition}
\newtheorem{definition&proposition}[theorem]{Definition and Proposition}
\newtheorem{definition&theorem}[theorem]{Definition and Theorem}
\newtheorem{example}[theorem]{Example}
\newtheorem{examples}[theorem]{Examples}
\newtheorem{definition}[theorem]{Definition}
\newtheorem{remarks}[theorem]{Remarks}
\newtheorem{remark}[theorem]{Remark}
\newtheorem{remarks*}[theorem]{Remarks$^*$}
\newtheorem{exercise}[theorem]{Exercise}
\newtheorem{exercise*}[theorem]{Exercise$^*$}
\newtheorem{exercises}[theorem]{Exercises}
\newtheorem{stelling}[theorem]{Stelling}
\newtheorem{teorema}[theorem]{Teorema}
\newtheorem{definisie}[theorem]{Definisie}
\newtheorem{voorbeelde}[theorem]{Voorbeelde}
\newtheorem{voorbeeld}[theorem]{Voorbeeld}
\newtheorem{opmerkinge}[theorem]{Opmerkinge}
\newtheorem{opmerking}[theorem]{Opmerking}
\newtheorem{oefeninge}[theorem]{Oefeninge}
\newtheorem{oefening}[theorem]{Oefening}

\def\<{\langle}

\def\>{\rangle}
\def\bproof{\noindent{\bf Proof: }}
\def\eproof{\begin{flushright}$\dashv$\end{flushright}}
\def\endbox{\begin{flushright}$\Box$\end{flushright}}
\def\beq{\begin{equation}}
\def\eeq{\notag \end{equation}}
\def\barR{\bar{\mathbb R}}
\def\=d{\buildrel\text{d}\over=}
\def\cadlag{c\`adl\`ag }


\fancyhead[LE,RO]{\thepage}
\fancyhead[LO,RE]{\em The L\'evy--It\^o Decomposition Theorem}
\section{L\'evy processes}
\begin{definition}\label{def_PAI}\rm Let $(\Omega,\mathcal F,\mathbb P, (\mathcal F_t)_t)$ be a filtered probability space. A stochastic process $X=(X_t)_t$ with values in $\mathbb R^n$ is said to be an $n$-dimensional L\'evy process if
\begin{enumerate}[(a)]\item $X$ is adapted to $(\mathcal F_t)_t$ (i.e. $X_t$ is $\mathcal F_t$--measurable for all $t\geq 0$).
\item $X_0=0$ a.s.
\item $X_{t+s}-X_t$ is independent of $\mathcal F_t$, and has the same law as $X_s$, for all $s,t\geq 0$.
\item $(X_t)_t$ is continuous in probability.
\end{enumerate}\endbox
\end{definition}

{\bf Consequences:}  Put $\varphi_t(u):=\mathbb E[e^{i\<u,X_t\>}]$, where $\<u,v\>$ represents the inner product on $\mathbb R^n$. From (d), the $\varphi_t(u)$ are continuous in the pair $(t,u)$.  From (b),  $\varphi_0(u)=1$. From (c), $\varphi_{t+s}(u)=\varphi_{t}(u)\varphi_s(u)$, so $\varphi_t(u)\not=0$ for any $(t,u)$. One can thus write $\varphi_t(u)=e^{-t\psi(u)}$, where $\psi$ is a continuous function null at 0.

Conversely, if $\psi$ is a continuous function null at 0, and such that for all $t\geq 0$, $\varphi_t(u):=e^{-t\psi(u)}$ is {\em positive definite} (which means that for each choice of finitely many $u_j,\lambda _j$ we have $\sum_{i,j}\lambda_i\bar{\lambda_j}\varphi_t(u_i-u_j)\geq 0$\footnote{I.e. for all ${u_1,\dots, u_k}$ the $k\times k$--matrix $(\varphi_t(u_i-u_j))_{i,j}$ is non--negative definite.}), then, by Bochner's theorem, $\varphi_t(u)$ is the Fourier transform  of a probability measure on $\mathbb R^n$. One can therefore construct a projective family of measures on $(\mathbb R^n)^{\mathbb R^+}$ by the formula
\begin{align*}
\mathbb E[e^{i\<u,X_{t_1}\>+\dots+i\<u,X_{t_n}\>}]&= \mathbb E[e^{i\<u_1+\dots+u_n, X_{t_1}\>+\dots+i\<u_n, X_{t_n}-X_{t_{n-1}}\>}]\\
&=\varphi_{t_1}(u_1+\dots+u_n)\cdot\varphi_{t_2-t_1}(u_2+\dots+u_n)\dots \varphi_{t_n-t_{n-1}}(u_n)\end{align*} for all finite choices of $0\leq t_1\leq\dots\leq t_n$. The process $(X_t)_t$ on $(\mathbb R^n)^{\mathbb R^+}$, given by Kolmogorov's Theorem, clearly adapted to the natural filtration $\mathcal F_t:=\sigma(X_s:s\leq t)$, possesses the properties (a), (b), (c), (d),  because the fact that $X_{t+s}-X_t\to 0$ in probability when $s\downarrow 0$ is an immediate  consequence\footnote{ Because then $X_s\buildrel{d}\over\to 0$, and  sequence of random variables which converges in distribution to a constant also converges in probability to that constant.} of the fact that $\varphi_s(u)\to 1$ when $s\downarrow 0$. Thus there is a L\'evy process to each $\psi$ possessing those properties.

\begin{theorem}\label{theorem_1} Suppose that $S$ is a countable subset of $\mathbb R^+$. Then there exists a $\mathbb P$--null set $N$ such that on $N^c$ the map $t\mapsto X_t$ has left- and right limits along $S$ (l\`ag, l\`ad). If one defines $Y_t:= \lim\limits_{s\in S, s\downarrow t} X_s$ on $N^c$, and 0 on $N$, then $Y$ is adapted to $\bar{\mathcal F}_t$, where $\bar{\mathcal F}_t$ is the completion of $\mathcal F_t$ in $\mathcal F$, i.e. completed by all the null sets {\em of $\mathcal F$} (or $\bigvee_{t}\mathcal F_t$).  Moreover, $Y$ is c\`adl\`ag ({\em continue \`a droite, pourvu limites \`a gauche}). Finally, $Y$ is a  modification of $X$, i.e. for all $t$ ,$\mathbb P(X_t\not= Y_t)=0$.
\end{theorem}
\bproof Suppose that $u\in\mathbb Q^n$, and that $M^u_t$ is defined by $M^u_t:=\frac{e^{i\<u,X_t\>}}{\varphi_t(u)}$. For each $u$, $(M^u_t)_t$ is a (complex) martingale, thus, save for a null set $N^u$, l\`agl\`ad along $S$ (see for example Neveu p. 129-132)\footnote{ Every martingale has a c\`adl\`ag modification.} Left and right limits along $S$ therefore  exist simultaneously for all the $M^u$, except for the null set $N:=\bigcup_{u\in\mathbb Q^n} N^u$. Suppose that, for $\omega\in N^c$, the map $s\mapsto X_s(\omega)$ could have two distinct cluster points $a, b$ as $s$ tends ($\uparrow$ or $\downarrow$) to $t$. One can always find a $u\in\mathbb Q^n$ such that $\<u,b-a\>\not\in 2i\pi\mathbb Z$; hence this is impossible\footnote{ $s\mapsto\frac{ e^{i\<u,X_s(\omega)\>}}{\varphi_s(u)}$ has unique left-- and right limits as $s\uparrow t$ or $s\downarrow t$ along $S$. If, e.g. we two have sequences $s_n\uparrow t$ and $s_n'\uparrow t$ such that $X_{s_n}(\omega\to a$ and $X_{s_n'}(\omega)\to b$, then $\frac{ e^{i\<u,X_{s_n}(\omega)\>}}{\varphi_{s_n}(u)}\to \frac{ e^{i\<u,a\>}}{\varphi_t(u)}$ and $\frac{ e^{i\<u,X_{s_n}(\omega)\>}}{\varphi_{s_n'}(u)}\to \frac{ e^{i\<u,b\>}}{\varphi_t(u)}$, so that $e^{i\<u,a\>}=e^{i\<u,b\>}$.} . Hence $X$ is l\`agl\`ad along $S$ on $N^c$, and $Y$ is c\`adl\`ag. By dominated convergence, we have  $\mathbb E[e^{i\<u,Y_t-X_t\>}]=\lim_{s\downarrow t}\mathbb E[e^{i\<u, X_s-X_t\>}]=1$, because $X_t$ is continuous in probability. Thus $\mathbb P(Y_t\not= X_t)=0$\footnote{ If $\mathbb E[e^{i\<u,X\>}]=1$ for all $u$, then $X$ has law $\delta_0$ by L\'evy inversion. So $\mathbb P(X=0)=1$.}. Finally, Because $X_t=Y_t$ a.s., $Y$ is adapted to $(\bar{\mathcal F}_t)_t$.\eproof

{\bf Consequence:} If we have a L\'evy process in the sense of Definition \ref{def_PAI}, one can now take for $X$ the regularised (i.e. c\`adl\`ag) version $Y$, and for $\mathcal F_t$ one can take $\overline{\sigma(X_s:s\leq t)}$. We now study a fixed L\'evy process $X$ (if one exists), and also fix $(\mathcal F_t)_t$ for the remainder of this chapter. 

\begin{theorem}\label{theorem_2} {\rm (0-1 Law)}\\
If we define $\mathcal F_{t^+}:=\bigcap_{s>t}\mathcal F_s$, then $\mathcal F_{t^+}=\mathcal F_t$
\end{theorem}

\bproof $\mathcal F_{t+}$ may be considered as a countable intersection, since $\mathcal F_{u}\subseteq\mathcal F_v$ when $u\leq v$. Thus if $t_1\leq t_2$, we have $\mathbb E[e^{i\<u,X_{t_1}\>}|\mathcal F_{t_2}] = \mathbb E[e^{i\<u,X_{t_1}\>}|\mathcal F_{t_2^+}] $, a common version being $e^{i\<u,X_{t_1}\>}$.

If $t_1>t_2$, then\footnote{ We use here the L\'evy-Doob Downward Theorem: If $Z$ is an integrable R.V. and $\mathcal G_n\downarrow\mathcal G$, then $\mathbb E[Z|\mathcal G_n]\to \mathbb E[Z|\mathcal G]$ a.s. and in $L^1$.} \begin{align*}\mathbb E[e^{i\<u,X_{t_1}\>}|\mathcal F_{t_2^+}] &\buildrel{\text{a.s.}}\over=\lim\limits_{s\downarrow t_2}\mathbb E[e^{i\<u,X_{t_1}\>}|\mathcal F_{s}] \\&
\buildrel{\text{a.s.}}\over=\lim\limits_{s\downarrow t_2}e^{i\<u,X_s\>}\varphi_{t_1-s}(u)\\&\buildrel{\text{a.s.}}\over=e^{i\<u,X_{t_2}\>}\varphi_{t_1-t_2}(u)\\&\buildrel{\text{a.s.}}\over=\mathbb E[e^{i\<u,X_{t_1}\>}|\mathcal F_{t_2}]
 \end{align*}
Thus for all $u,s$ we have $\mathbb E[e^{i\<u,X_s\>}|\mathcal F_{t+}]\buildrel{\text{a.s.}}\over=\mathbb E[e^{i\<u,X_s\>}|\mathcal F_t]$.  The two conditional expectations are equal a.s.  for all random variables $e^{i\<u,X_s\>}$, and hence for all\footnote{ We see that $\mathbb E[e^{i\sum_{j}\<u_j X_{s_j}\>}|\mathcal F_{t}]= \mathbb E[e^{i\sum_{j}\<u_j X_{s_j}\>}|\mathcal F_{t^+}]$  a.s. for all $u_j,s_j$. Now the collection $\mathcal M:=\{e^{i\sum_{j}\<u_j X_{s_j}\>}:u_j,s_j\}$ is closed under multiplication and conjugation. If $\mathcal H:=\{Z: \mathbb E[Z|\mathcal F_{t^+}]=\mathbb E[Z|\mathcal F_t]\}$, then $\mathcal H$ is a vector space closed under bounded pointwise convergence. Hence $\mathcal H$ contains every $\sigma(\mathcal M)$--measurable function  --- this is the monotone class theorem for complex-valued functions --- cf. Bichteler {\em Stochastic Integration with Jumps}, Exercise A.3.5.}
$\bigvee_t\mathcal F_t$--measurable random variables.\eproof

Henceforth, therefore, $(\mathcal F_t)_t$ is right--continuous: $\mathcal F_{t^+}=\mathcal F_t$. (In particular, $A\in\mathcal F_0=\mathcal F_{0^+}$ implies $\mathbb P(A) = 0$ or 1.)

\begin{theorem}\label{theorem_3} {\rm (Strong Markov property)}\newline If $T$ is a stopping time, then on $\{T<\infty\}$ the process $(X_{T+t}-X_T)_{t\geq 0}$ is a L\'evy process with the same law as $X$, adapted to $(\mathcal F_{T+t})_t$, c\`adl\`ag, and independent of $\mathcal F_T$.
\end{theorem}

\bproof
Suppose first that $T$ is bounded, and let $A\in\mathcal F_T$. Take $u_j\in\mathbb Q^n$, and $t\in\mathbb R^+$. Then 
\[\mathbb E\left[I_Ae^{i\sum_j\<u_j,X_{T+t_j}-X_{T+t_{j-1}}\>}\right] =\mathbb P(A)\prod_j\varphi_{t_j-t_{j-1}}(u_j)\] on application of the optional sampling theorem  to the martingales $M_t^{u_j}$.\footnote{ \begin{align*}\mathbb E\left[I_Ae^{i\sum_{j=1}^m\<u_j,X_{T+t_j}-X_{T+t_{j-1}}\>}\right]&= \mathbb E\left[I_A\prod_{j=1}^m\frac{M^{u_j}_{T+t_j}}{M^{u_j}_{T+t_{j-1}}}\frac{\varphi_{T+t_j}}{\varphi_{T+t_{j-1}}}\right] \\&=\mathbb E\left[I_A\prod_{j=1}^{m-1}\frac{M^{u_j}_{T+t_j}}{M^{u_j}_{T+t_{j-1}}}\frac{\varphi_{T+t_j}}{\varphi_{T+t_{j-1}}}\mathbb E\left[\frac{M^{u_m}_{T+t_m}}{M^{u_m}_{T+t_{m-1}}}\frac{\varphi_{T+t_m}}{\varphi_{T+t_{m-1}}}|\mathcal F_{T+t_{m-1}}\right]\right] \end{align*} Now  $\frac{\varphi_{T+t_m}}{\varphi_{T+t_{m-1}}}=\varphi_{t_m-t_{m-1}}(u_m)$ and $\mathbb E\left[\frac{M^{u_m}_{T+t_m}}{M^{u_m}_{T+t_{m-1}}}|\mathcal F_{T+t_{m-1}}\right]= 1$. Apply $m$ times.} If $T$ is  not bounded, the formula remains true when applied to $T\land n$ and $A\cap\{T\leq n\}\in\mathcal F_{T\land n}$. One can pass to the limit by dominated convergence, and hence the formula is true without restrictions. On the one hand it shows that $X_{T+t}-X_T$ is independent of $\mathcal F_T$, and on the other that $X_{T+t}-X_T$ has properties (a), (b), (c). It is clear that $X_{T+t}-X_T$ is c\`adl\`ag, thus {\em a fortiori} continuous in probability.\footnote{ With $A=\Omega$, the formula  shows that the law of $(X_t)_t$ and $(X_{T+t}-X_T)_t$ are the same.}
\eproof

\begin{corollary}\label{corollary_3} A L\'evy process which has  amplitudes of discontinuities a.s. bounded has moments of all orders.
\end{corollary} 
\bproof
Suppose $M$ is such that $\mathbb P(\exists t:|X_t-X_{t-}|\geq M)=0$. Put $T_1:=\inf\{t||X_t|\geq M\}$, and $T_n:=\inf\{t: t>T_{n-1}, |X_t|\geq M\}$. The right-continuity implies that the $T_n$ form a strictly increasing sequence of stopping times. Since $|X_T-X_{T-}|\leq M$ for all $T$, by iduction we have $\sup\limits_{s\leq T_n}|X_s|\leq 2nM$, the strong Markov property implies that $T_n-T_{n-1}$ is independent of $\mathcal F_{T_{n-1}}$, with the same law as $T_1$, and hence $\mathbb E[e^{-T_n}]=\mathbb E[e^{-T_1}]^n=a^n$, where $a<1$. Thus $\mathbb P(|X_t|\geq 2nM\}\leq \mathbb P(T_n<t)\leq a^n$, and hence there exists an exponential moment for $X_t$.\footnote{Choose $b$ s.t. $0<b<-\frac{\ln a}{2M}$. Then $e^{2bM}a<1$, so $\sum_n (e^{2bM}a)^n<\infty$. Now \[\aligned \mathbb E[e^{b|X_t|}]&=\sum_n\mathbb E[e^{b|X_t|}\Big|2(n-1)M<|X_t|\leq 2nM]\mathbb P(2(n-1)M<|X_t|\leq 2nM)\\&\leq \sum_ne^{2nMb}a^{n-1}=\frac1a\sum_n(e^{2bM}a)^n<\infty\endaligned\] Finally, $\mathbb E[e^{b|X_t|}]\geq\frac{b^n|X_t|^n}{n!}$, so $\mathbb E[|X_t|^n]<\infty$.}
\eproof
\section{ Poisson Process}
This is an increasing adapted  L\'evy process  which grows only by jumps of amplitude 1. We will denote it by $(N_t)_t$ in what follows, with or without supplementary indices.

If $T_1:=\inf\{t:N_t\not=0\}$, then $\{T_1>t\} = \{N_t=0\}$. $T_1$ is a stopping time, $\mathbb P(T_1>t+s) =\mathbb P(N_{t+s}-N_t=0, N_t=0)$, so by the strong Markov property, $\mathbb P(T_1>t+s)=\mathbb P(T_1>s)\mathbb P(T_1>t)$. This function being decreasing and bounded, we  have $\mathbb P(T_1>t)=e^{-at}$ for some $a\in\mathbb R^+$ ($T_1>0$ a.s.).
For $a=0$, $N_t\equiv 0$; if not, $T_1$ is a.s. finite, and if we put $T_n-T_{n-1}:= \inf\{t>0: N_{t+T_n}-N_{t+T_{n-1}}>0\}$, then $T_n-T_{n-1}$ is independent of $\mathcal F_{T_{n-1}}$ and has the same law as $T_1$. 

Then $\mathbb P(N_t=n)=\mathbb P(T_{n+1}>t, T_n\leq t) =\frac{a^nt^n}{n!}e^{-at}$, \footnote{  By induction, $T_n$ has density function $t\mapsto \frac{a^nt^{n-1}}{(n-1)!}e^{-at}$, so \[\mathbb P(T_{n+1}>t, T_n\leq t)=\int_0^t\mathbb P(T_{n+1}-T_n|t-s|T_n=s)\frac{a^ns^{n-1}}{(n-1)!}e^{-as}\;ds =\int_0^t e^{-a(t-s)}\frac{a^ns^{n-1}}{(n-1)!}e^{-as}\;ds= \frac{a^nt^n}{n!}e^{-at}\]}, and  $\mathbb E[e^{iuN_t}]=e^{-at(1-e^{iu})}$. As this function is  positive definite\footnote{ $\sum_{j,k}\lambda_j\bar{\lambda _k}e^{-at(1-e^{i(u_j-u_k)})}=e^{-at}|\sum_j\lambda _je^{iu_j}|^2$.}, then by the {\em converse} on p. 1, there exists a Poisson process\footnote{ As on p.1, Kolmogorov's theorem guarantees the existence of a stochastic process with the correct law, and Theorem 1.2 guarantees the existence ofa c\`adl\`ag version thereof.} Finally (by Corollary 1.5), $\hat{N}_t:=N_t-at$ and $(N_t-at)^2-at$ are integrable, and are martingales, as one immediately verifies\footnote{\label{footnote_Poisson_mean} Using the characteristic function $\mathbb E[e^{iuN_t}]=e^{-at(1-e^{iu})}$, and differentiating w.r.t. $u$, it is easy to see  that $\mathbb E[N_t]=at=\text{Var}(N_t)$.}

\begin{theorem}\label{theorem_4} If $M$ is a centered square--integrable martingale, $N$ a Poisson process, then for all $t$,\[\mathbb E[M_tN_t]=\mathbb E\left[\sum_{n\geq 0}(M_{T_n}-M_{T_n-})I_{\{T_n\leq t\}}\right]\] (where $T_n$ are the jump times of $N$.)
\end{theorem}
\bproof Suppose that $0=t_0<t_1<t_2<\dots<t_n=t$ is a partition of $[0,t]$. By using the  martingale property of $M_t$ and $\hat{N_t}:=N_t-at$ repetitively, we obtain:
\begin{align*} \mathbb E[M_tN_t]&=\mathbb E[M_t\hat{N}_t]=\mathbb E\left [\sum_i(M_{t_{i+1}}-M_{t_i})\sum_j(\hat{N}_{t_{j+1}}-\hat{N}_{t_j})\right]
\\&= \mathbb E\left [\sum_i(M_{t_{i+1}}-M_{t_i})(\hat{N}_{t_{i+1}}-\hat{N}_{t_i})\right]= \mathbb E\left [\sum_i(M_{t_{i+1}}-M_{t_i})({N}_{t_{i+1}}-{N}_{t_i})\right]
\end{align*}
If the step size $\sup_{i}(t_{i+1}-t_i)$ tends to $0$, \[ \sum_i(M_{t_{i+1}}-M_{t_i})({N}_{t_{i+1}}-{N}_{t_i})\buildrel{\mathbb P\text{ or a.s. }}\over\longrightarrow \sum_{n\geq 0}(M_{T_n}-M_{T_n-})I_{\{T_n\leq t\}}\]
The proof is complete if one can show that the Lebesgue dominated convergence theorem is applicable: Now \[\left|\sum_i(M_{t_{i+1}}-M_{t_i})(\hat{N}_{t_{i+1}}-\hat{N}_{t_i})\right|\leq 2\sup_{s\leq t}|M_s|\cdot N_t\] and both factors are in $L^2$ ($\mathbb E[\sup_{s\leq t}|M_s|^2]\leq 4\mathbb E[|M_t|^2]$) \footnote{Doob's $L^2$--inequality. Moreover the product of two $L^2$--variables is in $L^1$, by H\"older's inequality.}.
\eproof

\section{The Decomposition of Paul L\'evy}
\subsection{Jump Measure}\label{section_jump_measure}
Let $B$ be a Borel set in $\mathbb R^n$ with $0\not\in\bar{B}$. By recursion, we define the stopping times
\[S^1_B:=\inf\{t>0: X_t-X_{t-}\in B\}\qquad S^n_B:=\inf\{t>S^{n-1}_B: X_t-X_{t-}\in B\}\]
One easily verifies that , because of right--continuity, $X(t,\omega)$ is jointly measurable in $(t,\omega)$, hence the $S^n_B$ are stopping times adapted to $\mathcal F_{t+}$, thus to $\mathcal F_t$, by the 0-1 law. \footnote{ It follows from the theory of capacities and analytic sets that if $\sigma$ is a stopping time, $B$ a Borel set and $Y$ a progressively measurable process, then $\tau:=\inf\{t>\sigma: Y_t\in B\}$ is a stopping time, provided that the filtration satisfies the usual conditions. Now if $X$ is c\`adl\`ag adapted, then $\Delta X:=X-X_{-}$ is progressively  measurable.} The right-continuity implies that $S^1_B>0$ a.s. and that $N_t(B):=\sum\limits_{n\geq 0}I_{\{S^n_B\leq t\}}<\infty$ a.s. (If not, there would be a discontinuity of the second kind on the trajectory $X_t$\footnote{ $f$ has a discontinuity of the second kind at $t$ if one of the limits $\lim\limits_{s\uparrow t}f(s)$, $\lim\limits_{s\downarrow t}f(s)$ does not exist. In this case, there would be infinitely many jumps of amplitude $>\varepsilon$ by time $t$, for some $\varepsilon >0$.}.) $N_t(B)$ is therefore a Poisson process (see later), and we denote  by $\nu(B)$ the parameter $\mathbb E[N_1(B)]$.\footnote{ If $N$ is a Poisson process with parameter $a$, then $\mathbb E[N_1]=a$; cf. footnote \ref{footnote_Poisson_mean}.} For each $\omega$, $N_t(dx,\omega)$ defines a $\sigma$--finite measure 
on $\mathbb R^n-\{0\}$, thus $\nu(dx):=\mathbb E[N_1(dx)]$ is equally a $\sigma$--finite measure $(\geq 0$) on $\mathbb R^n-\{0\}$.

\subsection{Associated Jump Processes}

\begin{lemma} \label{lemma_1}
Let $f$ be  bounded measurable on $B$ in $\mathbb R^p$. Then
\[\int_Bf(x)\;N_t(dx) =\sum_{n\geq 1} f(X_{S^n_B}-X_{S^n_B-})I_{\{S^n_B\leq t\}}\]
\end{lemma}

\bproof If $f$ is a step function , $f=\sum_ja_j I_{B_j}$ where $\sum_j I_{B_j}=I_B$, the integral is $\sum_j a_j N_t(B_j) = \sum_j a_j(\sum_n I_{\{S^n_{B_j}\leq t\}})$. But the family $\{S^n_B\}$ is the union of the $\{S^n_{B_j}\}$, and the result follows\footnote{Note that if  $\sum_jI_{B_j}=I_B$, then $f(\Delta X_{S^n_B(\omega)})I_{\{S^n_B(\omega)\leq t\}}=f(\Delta X_{S^m_{B_j}(\omega)})I_{\{S^m_{B_j}(\omega)\leq t\}}$ for some $j,m$, and conversely. Hence $\sum_{n\geq 1}f(\Delta X_{S^n_{B}})I_{\{S^n_B\leq t\}}=\sum_{n\geq 1}\sum_j f(\Delta X_{S^n_{B_j}})I_{\{S^{n}_{B_j}\leq t\}}$. Now if $f=\sum_ja_jI_{B_j}$, then $\int_Bf(x)\;N_t(dx) = \sum_ja_jN_t(B_j) = \sum_ja_j\sum_{n\geq 1} I_{\{S^n_{B_j}\leq t\}} =\sum_{n\geq 1}\sum_ja_jI_{\{S^n_{B_j}\leq t\}} =\sum_{n\geq 1}\sum_j f(\Delta X_{S^n_{B_j}})I_{\{S^n_{B_j}\leq t\}}= \sum_{n\geq 1}f(\Delta X_{S^n_{B}})I_{\{S^n_{B}\leq t\}}$.} for step functions $f$. Else, one approximates $f$ uniformly by step functions\dots\footnote{If $f_k\to f$ uniformly, where the $f_k$ are step functions, then $\int_B f_k(x)\;N_t(dx)\to \int_B f(x)\;N_t(dx)$ a.s., by dominated convergence. Hence $\sum_{n\geq1} f_k(\Delta X_{S^n_B})I_{\{S^n_B\leq t\}}\to \int_Bf(x)\;N_t(dx)$ a.s. But for $t,\omega$, each sum\\ $\sum_{n\geq1} f_k(\Delta X_{S^n_B}(\omega))I_{\{S^n_B(\omega)\leq t\}}$ has a.s. only finitely many terms, so \\$\sum_{n\geq1} f_k(\Delta X_{S^n_B}(\omega))I_{\{S^n_B(\omega)\leq t\}}\to \sum_{n\geq 1}f(\Delta X_{S^n_B}(\omega))I_{\{S^n_B(\omega)\leq t\}}$.} \eproof

\begin{remark}\rm
In fact, for $B$ a Borel set ($0\not\in\bar{B}$ of course), it suffices that $f$ be finite everywhere on $B$ for the formula to be true, for $N_t(B)$ is finite a.s. for all $t$. In particular, we denote by $X_t(B)$ the quantity
\[\int_Nx\;N_t(dx) = \sum_{n\geq 1} (X_{S^n_B}- X_{S^n_B-})I_{\{S^n_B\leq t\}}\]
\endbox
\end{remark}

\begin{lemma}  \label{lemma_2}$\int_Bf(x)\;N_t(dx)$, $X_t(B)$ are L\'evy processes adapted to $(\mathcal F_t)_t$.
\end{lemma}
\bproof $N_t(dx)$ is an adapted L\'evy process!\eproof

\begin{lemma}  \label{lemma_3}$X_t-X_t(B)$ is a L\'evy process adapted to $(\mathcal F_t)_t$.\endbox
\end{lemma}

[To demonstrate that $N_t(B), X_t(B)$ and $X_t-X_t(B)$ are adapted L\'evy processes, we note that conditions (a) and (b) are automatically satisfied, and that (d) follows from the c\`adl\`ag property\footnote{Because $X$ has stationary increments, we have $\lim\limits_{s\uparrow t}\mathbb P(|X_t-X_s|>\varepsilon)=\lim\limits_{u\downarrow 0}\mathbb P(|X_u|>\varepsilon)=\lim\limits_{s\downarrow t}\mathbb P(|X_t-X_s|>\varepsilon)$. Since $X$ is c\`adl\`ag, $X_u\buildrel{\text{a.s.}}\over\to 0$, as $u\downarrow 0$, and hence also $X_u\buildrel{\mathbb P}\over\to 0$.}.
Only (c) remains to be verified. Now let $Z_t$ be one of the above-mentioned processes. Note that $Z_{t+s}-Z_t\in\sigma(X_u:u\leq t\leq t+s)$ is independent of $\mathcal F_t$. The same goes for the stationarity of the increments\footnote{The process $\hat{Z}$ defined by $\hat{Z}_s:= Z_{t+s}-Z_t$ has the same law as  $Z$, by the Strong Markov property. E.g. $\hat{X}_s(B)$ is number of jumps $X$ has in $B$ between times $t$ and $s+t$, and this has the same law as $X_s(B)$.}\dots]

$X_t-\int_{\{|x|\geq 1\}}N_t(dx)$  does not have jumps of $|\text{amplitude}|\geq 1$ (by Lemma \ref{lemma_1}), is a L\'evy process adapted to $(\mathcal F_t)_t$ (by Lemma  \ref{lemma_3}), and can thus be centered by a translation $\gamma t$ (by Corollary \ref{corollary_3}). We may therefore restrict  ourselves to the study of:

\subsection{ The L\'evy Decomposition for centered L\'evy process with jumps bounded by 1}
\begin{lemma}\label{lemma_4} Suppose that  $B\subseteq \{|x|\geq 1\}$ is such that $0\not\in \bar{B}$, and that $f:\mathbb R^n\to\mathbb R$ is such that $fI_B$ is in $L^2(\nu(dx))$ (the jump measure $\nu(dx)$ has been introduced in \S\ref{section_jump_measure}). Then we have
\[\mathbb E\left[\int_Bf(x)\;N_t(dx)\right]=t\int_Bf(x)\;\nu(dx)\]
and\[\mathbb E\left[\left(\int_Bf(x)\;N_t(dx)-t\int_Bf(x)\;\nu(dx)\right)^2\right]= t\int_Bf(x)^2\;\nu(dx)\]
\end{lemma}
\bproof If $f$ is a step function, $f=\sum_ja_jI_{B_j}$, we have\footnote{Because $\nu(B):=\mathbb E[N_1(B)]$, it follows that $\mathbb E[N_t(B)]=t\nu(B)$: By stationary independent increments it is clear that for natural numbers $p,q$ we have $\mathbb E[N_p(B)]=p\nu(B)$, and that $\mathbb E[N_p(B)] = q\mathbb E[N_{\frac{p}{q}}(B)]$, so that $\mathbb E[N_{\frac{p}{q}}(B)]=\frac{p}{q}\nu(B)$ for any non-negative rational $\frac{p}{q}$. Now by the c\`adl\`ag property, $N_t(B)=\lim\limits_{r\in\mathbb Q, r\downarrow t} N_r(B)$, so dominated convergence yields $\mathbb E[N_t(B)]=t\nu(B)$.}
\[\mathbb E\left[\sum_ja_jN_t(B_j)\right]=\sum_ja_j \mathbb E[N_t(B_j)]=t\sum_j a_j\nu(B_j)\]
For the second equation, note that if $B_i\cap B_j=\varnothing$, then by Theorem \ref{theorem_4} we have $\mathbb E[\hat{N}_t(B_i)\hat{N}_t(B_j)]=0$\footnote{Observe that $\mathbb E\left[\left(\sum_ja_j(N_t(B_j)-t\nu(B_j)\right)^2\right]=\mathbb E\left[\left(\sum_j\hat{N}_t(B_j)\right)^2\right]=\sum_ja_j^2\mathbb E[\hat{N}_t(B_j)^2]=\sum_ja_j^2t\nu(B_j)$, using footnote \ref{footnote_Poisson_mean}.}. For $f$ not a step function, choose a sequence of step functions $f_n$ such that $f_nI_B$ tends to $fI_B$ in $L^2(d\nu)$, and thus also in $L^1(d\nu)$. We then have  convergence of the corresponding stochastic integrals in $L^2(d\mathbb P)$ and $L^1(d\mathbb P)$.\footnote{Choose a sequence $f_n$ of step functions so that $|f_n-f|\leq 2^{-n}$. Put $Z:=\int_Bf(x)\;N_t(dx)$, $Z_n:=\int_Bf_n(x)\;N_t(dx)$. By dominated convergence, we have $|Z_n(\omega)-Z(\omega)|\leq \int_B|f_n(x)-f(x)|\;N_t(dx,\omega)\leq 2^{-n}N_t(B)\to 0$, so by dominated convergence $\mathbb E[Z_n]\to\mathbb E[Z]$. Similarly, $\int_Bf_n(x)\;\nu(dx)\to \int f(x)\;\nu(dx)$. Since $\mathbb E[Z_n]=t\int_Bf_n(x)\;\nu(dx)\to t\int_Bf(x)\;\nu(dx)$, we obtain $\mathbb E[Z]=t\int_Bf(x)\;\nu(dx)$.  Similarly, let $Y:=\int_Bf(x)\;N_t(dx)-t\int_Bf(x)\;\nu(dx)$, $Y_n:=\int_Bf_n(x)\;N_t(dx)-t\int_Bf_n(x)\;\nu(dx)$. Then $|Y_n-Y|\leq 2^{-n}(N_t(B)+t\nu(B))$, so by dominated convergence $\mathbb E[|Y_n-Y|^2]\to 0$, so $\mathbb E[Y_n^2]\to\mathbb E[Y^2]$. Similarly,
$\int_Bf_n(x)^2\;\nu(dx)\to \int_Bf(x)^2\;\nu(dx)$. Now $\mathbb E[Y^2]=t\int_Bf_n(x)^2\;\nu(dx)\to t\int_Bf(x)^2\;\nu(dx)$.}
\eproof

We now introduce $\mathcal M$, the space of c\`adl\`ag centered square-integrable martingales on $(\Omega,\mathcal F,\mathbb P)$, adapted to $(\mathcal F_t)_t$. We equip this space with the topology (of Fr\'echet) induced by the family of seminorms $q_t(M):=\mathbb E[M_t^2]^\frac12$. From the classical inequality\footnote{The Doob $L^2$--inequality.} $\mathbb E[\sup\limits_{s\leq t}M_s^2]\leq 4\mathbb E[M_t^2]$ we deduce that the $q_t$--convergence of  a sequence implies, with probability 1, the uniform convergence of the trajectories on the interval $[0,t]$, and thus the  limit is c\`adl\`ag. $q_t$-convergence equally implies   convergence of the random variables in $L^2$, and thus preserves the centeredness and martingale properties. In other words, $\mathcal M$ is closed in its topology.
\begin{lemma}\label{lemma_5} If $B$ is as above, and \[\mathcal H_B:=\left\{\int_Bf(x)\;N_t(dx)-t\int_Bf(x)\;\nu(dx):fI_B\in_2(d\nu)\right\}\] then $\mathcal H_B$ is a closed subspace of $\mathcal M$.
\end{lemma}
\bproof We have\footnote{A direct consequence of Lemma \ref{lemma_4}.}
\[t||fI_B||^2_{L^2(d\nu)}=q_t(M_{fI_B})^2\tag{$\star$}\]
where $M_{fI_B,t}:= \int_Bf(x)\;N_t(dx)-t\int_Bf(x)\;\nu(dx):fI_B\in_2(d\nu)$.
\newline $(\alpha)$: For $fI_B$ a step function, $M_{fI_B}$ is a martingale in $\mathcal M$, because to each Poisson process $N$ corresponds the martingale $\hat{N}_t:=N_t-\mathbb E[N_t]$ in $\mathcal M$. As all $L^2(d\nu)$--functions are limits of step functions, $M_{fI_B}$ is a martingale in $\mathcal M$ as soon as $fI_B\in L^2(d\nu)$.
\newline$(\beta)$: Now $\mathcal H_B$ is closed in $\mathcal M$, because $q_t$--convergence also implies convergence in $L^2(d\nu)$, by $(\star$)\footnote{Suppose that $M_n:=M_{f_nI_B}$ for some sequence $f_n\in L^2(d\nu)$, and that $(M_n)_n$ is a Cauchy sequence in $\mathcal H_B$. Then by $(\star)$, $(f_nI_B)_n$ is Cauchy in $L^2(d\nu)$ and thus converges to some $f=fI_B$. By $(\star)$ again, $M_n\to M_{fI_B}$ in $\mathcal H_B$.}
\eproof

\begin{lemma}\label{lemma_6} Let $B$ be as in Lemma \ref{lemma_4}. If $M\in\mathcal M$ is continuous at the times $S_B^n$, then $M$ is orthogonal to $\mathcal H_B$. 
\end{lemma}
\bproof By Theorem \ref{theorem_4}, for all $A\subseteq B$ and all $t$ we have $\mathbb E[M_tN_t(A)]=0$. Now, the $\{N_t(A)|A\subseteq B\}$ generate $\mathcal H_B$\footnote{The compensated Poisson processes $(\hat{N}_t(A))_t$, where $A$ is a Borel subset of $B$, generate $\mathcal H_B$: Each $fI_B\in L^2(d\nu)$ is a limit of step functions $\sum_ja_jI_{A_j}$, where $A_j\subseteq B$. Then $\int_Bf\;N_t(dx)-t\int_Bf(x)\;\nu(dx)$ is a limit of $\sum_ja_j\hat{N}_t(A_j)$.}.

\begin{corollary}\label{corollary_6} If $B_1$, $B_2$ are disjoint Borel sets, with $0\not\in\bar{B_1}\cup\bar{B_2}$, then the processes $(X_t(B_1))_t$ and $(X_t(B_2))_t$ are independent L\'evy processes.
\end{corollary}
\bproof That they are L\'evy processes has already been demonstrated (Lemma \ref{lemma_4}).
If now \[M^u_{1,t}:=\frac{e^{i\<u, X_t(B_1)\>}}{\mathbb E[e^{i\<u, X_t(B_1)\>}]}-1\qquad M^v_{2,t}:=\frac{e^{i\<v, X_t(B_2)\>}}{\mathbb E[e^{i\<v, X_t(B_2)\>}]}-1\] then these two martingales are orthogonal, by Lemma \ref{lemma_6}. We have $\forall s,t\in\mathbb R^+\;\forall u,v\in\mathbb R^n\;(\mathbb E[M^u_{1,t}M^v_{2,s}]=0)$, which ensures\footnote{Recall a result of Ka\v c, which states that two random variables $Y,Z$ are independent iff $\mathbb E[e^{i\<u,Y\>+i\<v,Z\>}]=\mathbb E[e^{i\<u,Y\>}]\mathbb E[e^{i\<v,Z\>}]$ for all $u,v$.} their independence.
\eproof

Now put \[Y_t(B):=X_t(B)-\mathbb E[X_t(B)]=X_t(B)-t\int_Bx\;\nu(dx)\]
Then $Y$ is  both a L\'evy process and a martingale in $\mathcal M$. If we define
\[B_k:=\left\{\frac{1}{k+1}<|x|\leq \frac{1}{k}\right\}\qquad\text{and}\qquad A_n:=\bigcup_{k=1}^nB_k\] then the $Y(B_k)$ are pairwise independent, and $X-Y(A_n)$ and $Y(A_n)$ are orthogonal, and even independent (retrace the proof of Corollary \ref{corollary_6}). Consequently, the series (sum) of the $Y(B_k)$ converges in $L^2$, and thus in $\mathcal M$, to a L\'evy process $X^{d}$, while $X-Y(A_n)$ converges to a L\'evy process  $X^{c}$ in $\mathcal M$.\footnote{Using independence, $q_t(\sum_{m<k\leq n}Y_t(B_k))=\sum_{m<k\leq n}\mathbb E[Y_t(B_k)^2]=t\int_{\bigcup_{m<k\leq n}B_k}|x|^2\;\nu(dx)\leq \frac{t}{m^2}\nu(\bigcup_{m<k\leq n}B_k)\to 0$ as $m\to \infty$. Hence the (sequence $(Y(A_n))_n$ corresponding to the) series $\sum_{k=1}^\infty Y(B_k)$ is Cauchy in $\mathcal M$, and thus converges to some $X^d\in\mathcal M$. Moreover, using dominated convergence and the fact that the $Y(A_n)$ are L\'evy processes, we obtain (i): $\mathbb E[e^{i\<u,X^d_{t+s}-X^d_t\>}|\mathcal F_t]=\lim_n\mathbb E[e^{i\<u,Y_{t+s}(A_n)-Y_t(A_n)\>}|\mathcal F_t]=\lim_n\mathbb E[e^{i\<u,Y_{t+s}(A_n)-Y_t(A_n)\>}]=
\mathbb E[e^{i\<u,X^d_{t+s}-X^d_t\>}]$, which shows that $X^d$ has independent increments, and (ii): $\mathbb E[e^{i\<u,X^d_{t+s}-X^d_t\>}]=\lim_n\mathbb E[e^{i\<u,Y_{t+s}(A_n)-Y_t(A_n)\>}]=\lim_n\mathbb E[e^{i\<u,Y_s(A_n)\>}]=\mathbb E[e^{i\<u,X^d_s\>}]$, which shows that $X^d$ has independent increments.}

Thus\footnote{Note that $Y(A_n)=\int_{\{\frac{1}{n+1}<|x|\leq 1\}} |x|\;(N_t(dx)-t\nu(dx))$, and that $X^d=\lim_nY(A_n)$}:
\begin{lemma}\label{lemma_7} $X_t=X^c_t+X^d_t$, where $X^c$ is a martingale with continuous sample paths, and \[X^d_t:=\int_{|x|\leq 1}x \;\left(N_t(dx)-t\nu(dx)\right)\]\endbox
\end{lemma}
\begin{remark}\rm This last integral  exists in $L^2$, and hence  $\int_{\{|x|\leq 1\}}|x|^2\;\nu(dx)<\infty$.
\endbox
\end{remark}
It remains to characterize the continuous part: We will show that it is necessarily Gaussian L\'evy process, i.e. that each $X^c_t$ is Gaussian. For this, it suffices to show this for each one--dimensional projection (a well-known property of the Gaussians\footnote{We start with three observations. (i): First observe that if $(X_t)_t$ is an $n$--dimensional L\'evy process and $A$ is a $m\times n$--matrix, then $AX$ is an $m$-dimensional L\'evy process: For $\mathbb E[e^{i\<u,AX_{t+s}-AX_t\>}|\mathcal F_t] =\mathbb E[e^{i\<A^\top u, X_{t+s}-X_t\>}|\mathcal F_t]=\mathbb E[e^{i\<A^\top u, X_s\>}]=\mathbb E[e^{i\<u,AX_s\>}]$ shows that $AX$ has independent stationary increments. (ii): Next recall that an $n$--dimensional  random vector $Y$ is multivariate Gaussian if and only if each linear combination $\lambda^\top Y$ is univariate Gaussian, for any $\lambda\in\mathbb R^n$. (iii): Finally, if  $(Z_t)_t$ is a n--dimensional L\'evy process such that each random vector $Z$ is multivariate Gaussian, then $(Z_t)_t$ is a Gaussian process: For each linear combination $\sum_{j=1}^m\lambda_kZ_{t_k}$ can be written as a linear combination of independent normally distributed random vectors  $\sum_{j=1}^m\lambda_kZ_{t_k}=\sum_{k=1}^m \gamma_k(Z_{t_k}-Z_{t_{k-1}})$ (where $\gamma_k:=\sum_{j=k}^m\lambda_j$ and $t_0:=0$), so that each linear combination $\sum_{j=1}^m\lambda_kZ_{t_k}$ is Gaussian, i.e. $(Z_t)_t$ is a Gaussian process. Having made these three observations, we now proceed: Suppose that $X=(X_t)_t$ is an $n$--dimensional L\'evy process with continuous sample paths. Then each linear combination $\lambda^\top X$ is a one--dimensional L\'evy process (by (i)) with continuous sample paths.  If we can show that, for each $t$ and $\lambda$, the random variable $\lambda^\top X_t$ is Gaussian, then by (ii) each random vector $X_t$ is multivariate Gaussian, so by (iii) $(X_t)_t$ is a Gaussian process. Thus it suffices to show that whenever $B$ is a one--dimensional L\'evy process with continuous sample paths, then each random variable $B_t$ is Gaussian ---  just take $B_t=\lambda^\top X_t$. }). In other words, we must show that\footnote{Here is another proof of Lemma \ref{lemma_8}, which uses the L\'evy characterization of Brownian motion. Suppose that $X_t$ is a $d$--dimensional L\'evy process with
continuous sample paths. Then it has moments of all orders, by
Corollary \ref{corollary_3}. In particular, the
process $X_t-\mathbb E[X_t]$ is a continuous martingale centered
at 0. We now show that any centered continuous L\'evy process is a
Brownian motion in the loose sense: Components need not be
independent, but are multi--variate Gaussian.

So let $X_t=(X^{(1)}_t,\dots, X^{(d)}_t)$ be a centered L\'evy
process with continuous sample paths. Let $A$ be the non--negative
definite symmetric $d\times d$--matrix defined by
\[A_{ij} = \mathbb E[X^{(i)}_1X^{(j)}_1]\] where $X_t =
(X^{(1)}_t,\dots, X^{(d)}_t)$. We claim that the quadratic covariation
process of $X^{(i)}_t$ and $X^{(j)}_t$ is given by
\[[X^{(i)},X^{(j)}]_t = A_{ij}t\]
Recall that $\mathbb E[e^{i\<u,X_t\>}]=e^{-t\psi(u)}$ for some
$\psi:\mathbb R^d\to\mathbb C$ with $\psi(0) = 0$.  Since $X_t$
has moments of all orders, the function $\psi$ is ${\cal
C}^\infty$.
Moreover, since $\mathbb E[X^{(i)}_t]=0$, we have $\frac{\partial
}{\partial u^{(i)}}\Big|_{u=0}e^{-t\psi(u)} = 0$, and thus that
$\frac{\partial }{\partial u^{(i)}}\Big|_{u=0}\psi(u) = 0$ also.
It now follows easily that
\[\mathbb E[X^{(i)}_tX^{(j)}_t] = -\frac{\partial^2}{\partial
u^{(i)}\partial u^{(j)}}\Big|_{u=0}e^{-t\psi(u)}
=-t\;\frac{\partial^2}{\partial u^{(i)}\partial
u^{(j)}}\Big|_{u=0}\psi(u)=t\;\mathbb E[X^{(i)}_1X^{(j)}_1]\]i.e.
that\[\mathbb E[X^{(i)}_tX^{(j)}_t]=A_{ij}t\] To show that
$[X^{(i)},X^{(j)}]_t = A_{ij}t$, it suffices to show that
$X^{(i)}_tX^{(j)}_t - A_{ij}t$ is a  martingale.
 By the fact that increments are independent with mean zero,
we have
\[\aligned {}&\phantom{llll}\mathbb
E[X^{(i)}_tX^{(j)}_t-X^{(i)}_sX^{(j)}_s|{\cal F}_s]\\&=\mathbb
E[(X^{(i)}_t-X^{(i)}_s)(X^{(j)}_t-X^{(j)}_s)+X^{(i)}_s(X^{(j)}_t-X^{(j)}_s)+X^{(j)}_s(X^{(i)}_t-X^{(i)}_s)|{\cal
F}_s]\\& = \mathbb
E[(X^{(i)}_t-X^{(i)}_s)(X^{(j)}_t-X^{(j)}_s)]\endaligned\] Taking
expectations on both sides shows that \[\mathbb
E[X^{(i)}_tX^{(j)}_t-X^{(i)}_sX^{(j)}_s] =\mathbb
E[(X^{(i)}_t-X^{(i)}_s)(X^{(j)}_t - X^{(j)}_s)]\] It follows that
\[\mathbb E[X^{(i)}_tX^{(j)}_t-X^{(i)}_sX^{(j)}_s|{\cal F}_s] = \mathbb
E[X^{(i)}_tX^{(j)}_t - X^{(i)}_sX^{(j)}_s] = A_{ij}(t-s)\] and
thus that $X^{(i)}_tX^{(j)}_t - A_{ij}t$ is a martingale.

Now, for $\lambda\in\mathbb R^d$, let $Z^{\lambda}_t = \<\lambda,
X_t\>$. Then $Z^{\lambda}_t$ is a centered continuous
one--dimensional martingale. Using the fact that the covariance
process bracket operation is bilinear, we see that the quadratic
variation of $Z^\lambda_t$ is given by
\[[Z^{\lambda}]_t = \<\lambda, A\lambda\>\;t\]
Hence, by L\'evy's characterization, $Z^{\lambda}_t$ is a Brownian
motion with variance parameter $\<\lambda, A\lambda\>$ (i.e.
$Z^{\lambda}_t \sim N(0, \<\lambda, A\lambda\>\;t)$). It now
follows that
\[\mathbb E[e^{i\<u,X_t\>}] =\varphi_{Z^u_t}(1)
=e^{-\frac12\<u,Au\>\;t}\] which proves that $X_t$ is a
$d$--dimensional Brownian motion with covariance matrix $A$.

Now if $X$ is a L\'evy process with continuous sample paths, then
$X_t-\mathbb E[X_t]$ is centered, and thus a Brownian motion. If
we define $\gamma=\mathbb E[X_1]$, then $\mathbb E[X_t] = \gamma
t$, and so $\mathbb E[e^{i\<u,X_t-\gamma t\>}] =
e^{-\frac12\<u,Au\>\;t}$ for some symmetric non--negative definite
matrix $A$. We have proved:
\begin{theorem} \label{thm_continuous_Levy_=_BM} Suppose that $X_t$ is a L\'evy process with
continuous sample paths. Then there exists $\gamma\in\mathbb R^d$
and a symmetric non--negative definite $d\times d$--matrix $A$
such that
\[\mathbb E[e^{i\<u,X_t\>}] =
e^{i\<u,\gamma\>\;t-\frac12\<u,Au\>\;t}\] Hence $X_t$ is a
$d$--dimensional arithmetic Brownian motion with drift $\gamma$
and covariance matrix $A$.\endbox\end{theorem}}:
\begin{lemma}\label{lemma_8} Let $B_t$ be a one--dimensional centered L\'evy process with continuous sample paths. Then there is $\sigma^2\in\mathbb R^+$ such that 
\[\mathbb E[e^{iuB_t}]=e^{-\frac12u^2\sigma^2t}\]
\end{lemma}

\bproof
A L\'evy process without discontinuities has moments of all orders, by Corollary \ref{corollary_3}. If $\mathbb E[B_t^2]=0$ for all $t>0$, then the problem is solved. If not, we can assume that $\mathbb E[B-t^2]=t$, by multiplying the process by a constant. Note that $\mathbb E[B^4_t]=at+bt^2+ct^3$: It suffices to put $\mathbb E[e^{iuB_t}]=e^{-t\psi(u)}$, to differentiate four times at the origin ($\psi(u)$ is of the class $\mathcal C^\infty$, like $\varphi_t(u)$) and to observe that $\psi'(0)=0$.\footnote{For $n\geq 1$, we have $\mathbb E[B_t^n]= i^{-n}\frac{d^n}{du^n}|_{u=0}e^{-t\psi(u)}=\sum_{k=1}^na_kt^k$. Here the $t^n$--term will have coefficient $a_n=(-\psi'(0))^n=0$.} Now let $P:=\{0=t_0<t_1<\dots<t_n=t\}$ be a partition of $[0,t]$ whose step size $\sup_{j}\{(t_{j+1}-t_j)$ will tend to 0. We denote by $\Delta t_j$ the quantity $t_{j+1}-t_j$, and by $\Delta B_j$ the quantity $B_{t_{j+1}}-B_{t_j}$. We then have
\begin{align*}
\mathbb E[e^{iuB_t}-1]&=\mathbb E\left[\sum_{j}e^{iuB_{t_{j+1}}}-e^{iuB_{t_j}}\right]\\&=
iu\sum_j\mathbb E\left[e^{iuB_{t_j}}\Delta B_j\right]-\tfrac12u^2\sum_j\mathbb E\left[e^{iuB_{t_j}}(\Delta B_j)^2\right]\\&-\tfrac12u^2\sum_j\mathbb E\left[(\Delta B_j)^2\cdot(e^{iu(B_{t_j}+\theta_j\Delta B_j)}-e^{iuB_{t_j}})\right]
\end{align*}where the $\theta_j$ are numbers between $0$ and $1$, by the second order Taylor formula. In the second line, the first term is zero: $\Delta B_j$ has zero expectation, and is independent of $B_{t_j}$.

The second term equates to $-\frac12u^2\sum_j \varphi_{t_j}(u)\Delta t_j$, and thus tends to $-\frac12u^2\int_0^t e^{-s\psi(u)}\;ds$ as the step size tends to 0. 

The third term tends to 0: Let $A_\alpha$ be the event \[A_\alpha:=\Big\{\sup\limits_j\sup\limits_{t_j\leq u,v\leq t_{j+1}} |B_u-B_v|<\alpha\Big\}\] The third term can therefore be bounded by\footnote{Observe that
\begin{align*}
\left|\tfrac12u^2\sum_j\mathbb E\left[(\Delta B_j)^2\cdot(e^{iu(B_{t_j}+\theta_j\Delta B_j)}-e^{iuB_{t_j}})\right]\right|&\leq\tfrac12|u|^2\int\sum_j\Delta B_j^2|(e^{iu(B_{t_j}+\theta_j\Delta B_j)}-e^{iuB_{t_j}}|\;d\mathbb P\\&\leq
\tfrac12|u|^2\int_{A_\alpha}\sum_j\Delta B_j^2|u\theta_j\Delta B_j|\;d\mathbb P +\tfrac12|u|^2\int_{A_\alpha^c}\sum_j\Delta B_j^2\cdot 2\;d\mathbb  P
\end{align*} using the facts that $|e^{ih}|\leq |h|$ and that $|e^{ih}-1|\leq 2$. Now $|u\theta_j\Delta B_j|\leq \alpha|u|$ on $A_\alpha$ }
\[\tfrac12|u|^3\int_{A_\alpha}\alpha (\sum_j\Delta B_j^2)\;d\mathbb P+|u|^2\int_{A_\alpha^c} \sum_j\Delta B_j^2\;d\mathbb P\]
thus\footnote{ $\int_{A_\alpha^c} \sum_j\Delta B_j^2\;d\mathbb P= \mathbb E[I_{A_\alpha^c}\sum_j\Delta B_j^2]\leq \sqrt{\mathbb E[I_{A_\alpha^c}^2]}\sqrt{\mathbb E[(\sum_j\Delta B_j^2)^2]}$, by the  Cauchy-Schwarz or H\"older inequality.} by
\[\tfrac12\alpha |u|^3\mathbb E[\sum_j\Delta B_j^2] +|u|^2\sqrt{\mathbb P(A_\alpha^c)}\cdot\sqrt{\mathbb E[(\sum_j\Delta B_j^2)^2]}\] thus, taking into account the evaluation of $\mathbb E[B_t^4]$, by $\tfrac12\alpha|u|^3t+u^2\sqrt{\mathbb P(A_\alpha^c)}(O(t+t^3))^\frac12$. Finally, note that that the continuity of the sample paths implies that, as the partition gets finer, $\mathbb P(A_\alpha^c)\to 0$. The expectation of the third term  has  limit $\leq \frac12\alpha |u|^3t$, and  thus converges to zero. We therefor obtain the equation:
\[e^{-t\psi(u)}-1 = -\tfrac12 u^2\int_0^te^{-s\psi(u)}\;ds\] which identifies $\psi(u)$ with $\frac12u^2$.\footnote{$-\tfrac12u^2\int_0^te^{-s\psi(u)}\;ds = \frac{u^2}{2\psi(u)}(e^{-t\psi(u)}-1)$.}
\eproof

\subsection{The Decomposition Theorem}
\begin{theorem}\begin{enumerate}[(A)]\item Let $X$ be an $n$--dimensional L\'evy process. Then \[X_t= B_t+t\mathbb E\left[X_1-\int_{|x|>1\}}xN_1(\;dx)\right]+\int_{\{|x|\geq 1\}}x\;N_t(dx) +\int_{\{|x|<1\}}x\Big(N_t(dx)-t\nu(dx)\Big)\] where \begin{itemize}\item $B_t$ is a centered Gaussian L\'evy  process with a.s. continuous sample paths.
\item $N_t(dx)$ is a family of Poisson processes, independent of $B_t$, with $N_t(A)$ idependent of $N_t(B)$ if  $A\cap B=\varnothing$, and with $\nu(dx)=\mathbb E[N_1(dx)]$.
\item $\nu(dx)$ is a positive measure on $\mathbb R^n-\{0\}$, with $\int |x|^2\land 1\;\nu(dx)<\infty$.
\item The first stochastic integral is in the sense of $L^0$ and the second in the sense $L^2$.
\end{itemize}
\item {\bf Formula for the law\footnote{This is the L\'evy--Khintchine formula for the characteristic function of an infinitely divisible distribution.}:} Under these conditions,
\begin{align*}\psi(u)&=-\tfrac{1}{t}\mathbb E[e^{i\<u,X_t\>}]\\&=\tfrac12 Q(u)-i\<a,u\>+\int_{\{|x|\geq 1\}} 1-e^{i\<u,x\>}\;\nu(dx)+\int_{|x|<1\}} 1-e^{i\<u,x\>} +i\<u,x\>\;\nu(dx)
\end{align*} where $Q$ is a positive definite quadratic form\footnote{i.e. $Q(u)=u^\top \Sigma u$ for some symmetric positive definite  matrix $\Sigma$. Here, $\Sigma$ is the covariance matrix for the continuous Gaussian part.} on $\mathbb R^n$, $a\in\mathbb R^n$, and $\nu(dx)$ is as in (A).
\item Conversely, given $Q,a,\nu$ as in (B), here exists a L\'evy process whose law is given by the formula in (B).
\item The representation in (B) is unique.
\end{enumerate}
\end{theorem}

\bproof To summarize the preceding, (B) is obvious. (D) The uniqueness of the decomposition is clear by constructions. For the law, put\[\psi(u):=\tfrac12  Q_j(u)-i\<a_j,u\>+\int_{\{|x|\geq 1\}}(1-e^{i\<u,x\>}\;\nu_j(dx) +\int_{\{|x|<1\}} 1-e^{i\<u,x\>}+i\<u,x\>\;\nu_j(dx)\] where $Q_j,a_j,\nu_j$ are as in (B), and $j=1,2$. Let $y$ be a unit vector in $\mathbb R^n$:x
$\lim\limits_{s\to \infty} \frac{\psi(sy)}{s^2} =\frac12Q_j(y)$, for $\lim\limits_{s\to\infty}\frac{\<a,sy\>}{s^2}=0$, and the same for the integral terms, by dominated convergence\footnote{\label{footnote_exponential_inequality}It is useful to observe that  for any $\theta\in\mathbb R$ and $n\in\mathbb N$, we have $|e^{i\theta}-\sum_{k=1}^{n-1}\frac{(i\theta)^k}{k!}|\leq \frac{|\theta|^n}{n!}$: This follows immediately from the identity $e^{i\theta}=\sum_{k=1}^{n-1}\frac{(i\theta)^k}{k!}+\frac{i^n}{(n-1)!}\int_0^\theta(\theta-x)^{n-1}e^{ix}\;dx$, which may be proved using induction and integration by parts. } \footnote{ Recall that $|x|^2\land 1$ in $\nu$--integrable. On $\{|x|\geq 1\}$ we have $\frac{1-e^{i\<sy,x\>}}{s^2}\leq \frac{2}{s^2}\leq 2$ for $s$ sufficiently large, and $2I_{\{|x|\geq 1\}}$ is $\nu$--integrable. On $\{|x|< 1\}$ we have $\frac{1-e^{i\<sy,x\>}-i\<sy,x\>}{s^2}\leq\frac{|\<sy,x\>|^2}{2s^2}\leq |x|^2$, and $|x|^2I_{\{|x|<1\}}$ is $\nu$--integrable.  In both cases, we have domination by a $\nu$--integrable function.  }. $\psi$ therefore determines $Q$. We will determine all the projections of $\nu$ (let us suppose that the L\'evy process is one--dimensional)\footnote{ For the one--dimensional case, we have \begin{align*} &\phantom{=}\psi(u)-\tfrac12\int_{u-1}^{u+1}\psi(v)-\tfrac12Q(v)\;dv\\
&=-iau +\int_{\mathbb R}1-e^{iux}-i uxI_{\{|x|<1\}}\nu(dx)  -\tfrac12\int_{u-1}^{u+1}-iav +\int_{\mathbb R}1-e^{ivx}-ivxI_{\{|x|<1\}}\nu(dx) \;dv\\
&=\tfrac12\left[\int_{-1}^1-iau +\int_{\mathbb R}1-e^{iux}-i uxI_{\{|x|<1\}}\nu(dx)\;dv -\int_{-1}^1-ia(u+v) +\int_{\mathbb R}1-e^{i(u+v)x}-i (u+v)xI_{\{|x|<1\}}\nu(dx)\;dv\right]\\
&=\tfrac12\left[\int_{-1}^1 aiv+\int_\mathbb R- (e^{iux}-e^{i(u+v)x})+(iux-i(u+v)x)I_{\{|x|<1\}}\;\nu(dx)\;dv\right]\\
&=-\tfrac12\left[\int_{-1}^1\int_\mathbb R e^{iux}(1-e^{ivx})+ivxI_{\{|x|<1\}}\;\nu(dx)\;dv\right]
\end{align*}
But, using the inequality derived in footnote \ref{footnote_exponential_inequality}, $| e^{iux}(1-e^{ivx})+ivxI_{\{|x|<1\}}|$ is $\leq 2$ on $\{|x|\geq 1\}$, and is $\leq |e^{iux}|\cdot|1-e^{ivx}+ivx|+|ivx|\cdot|1-e^{iux}|\leq \frac12|v|^2|x|^2+|v||x||ux|=K|x|^2$ on $\{|x|<1\}$. Hence $ e^{iux}(1-e^{ivx})+ivxI_{\{|x|<1\}}$ is $\nu(dx)\otimes dv$--integrable, so Fubini's Theorem applies:
\begin{align*}&\phantom{=} 
-\tfrac12\int_{-1}^1\int_\mathbb R e^{iux}(1-e^{ivx})+ivxI_{\{|x|<1\}}\;\nu(dx)\;dv\\&= -\tfrac12\int_\mathbb Re^{iux}\int_{-1}^1 1-e^{ivx}+ivxI_{\{|x|<1\}}\;dv\;\nu(dx)\\
&= -\int_\mathbb R e^{iux}\left(1-\tfrac{\sin x}{x}\right)\;\nu(dx)
\end{align*}
Hence $-\psi(u)+\tfrac12Q(u)$ is the Fourier transform of the positive measure $\rho(dx):=(1-\frac{\sin x}{x}\;\nu(dx)$. By Fourier inversion, $\rho$ is determined by $\psi$ and $Q$. Since $\psi$ determines $Q$, we see that $\rho$, hence $\nu$, is determined by $\psi$.}:
\[\psi(u)-\tfrac12Q(u)-\tfrac12\int_{u-1}^{u+1}\psi(v)-\tfrac12Q(v)\;dv = -\int e^{iux}(1-\tfrac{\sin x}{x})\;\nu(dx)\] The lefthand side therefore determines (Bochner) the (positive) measure $(1-\frac{\sin x}{x})\;\nu(dx)$, and thus $\nu$, because $(1-\frac{\sin x}{x})>0$ on $\mathbb R-\{0\}$. Now $a$ can be identified by subtraction.

(C): For each $t$, $e^{-tQ(u)}$, $e^{it\<a,u\>}$ and $e^{-t(1-e^{i\<u,x\>})}$ are positive definite\footnote{Each is clearly a characteristic function.}, thus the equation (B) defines\footnote{Again, using the inequality of footnote \ref{footnote_exponential_inequality}, we have $|1-e^{i\<u,x\>}|\leq 2$ on $\{|x|\geq 1\}$ and $\leq\frac12 |u|^2|x|^2$ on $\{|x|<1\}$, so if $\int |x|^2\land 1\;\nu(dx)$ is finite, both integrals in (B) are defined.} a continuous function if $\int |x|^2\land 1\;\nu(dx)<\infty$, with each $e^{-t\psi(u)}$ being positive definite, and thus, according to  the observations following the first paragraph, defines a L\'evy process.\eproof
\begin{remark}\rm We have shown, without studying this process, that the c\`ad\l`ag version of a Brwonain motion is a.s. continuous! In effect, $\psi(u)=-\tfrac12u^2$ defines a L\'evy process in in the sense of Definition \ref{def_PAI} because for each $t$, $e^{-t\frac{u^2}{2}}$ is positive definite. We regularize, remove the jumps, and obtain  a formula (B). Uniqueness shosw that $\nu(dx)\equiv 0$, $a=0$, and thus there are no discontinuities!\endbox
\end{remark}
\newpage
\theendnotes
\end{document}